\documentclass[12pt]{article}
\usepackage{amssymb}

\author{I. G. Korepanov\\
\normalsize Southern Ural State University\\[-0.5ex]
\normalsize 76 Lenin avenue\\[-0.5ex]
\normalsize 454080 Chelyabinsk, Russia\\[-0.5ex]
\normalsize E-mail: kig@susu.ac.ru}
\date{}
\title{Euclidean 4-simplices and invariants of four-dimensional
manifolds: II.~An algebraic complex and moves~$2\leftrightarrow 4$}
\def\be{\begin{equation}}
\def\ee{\end{equation}}

\def\ol{\overline}
\def\pa#1#2{{\partial#1\over\partial#2}}

\def\Im{\mathop{\rm Im}\nolimits}
\def\Ker{\mathop{\rm Ker}\nolimits}

\newtheorem{theorem}{Theorem}
\newtheorem{lemma}{Lemma} 

\def\emline#1#2#3#4#5#6{%
       \put(#1,#2){\special{em:moveto}}%
       \put(#4,#5){\special{em:lineto}}}
\def\newpic#1{}
\begin{document}
\maketitle

\begin{abstract}
%%%%%%%%%%%%%
%\baselineskip1.85\baselineskip
%%%%%%%%%%%%%
We write out some sequences of linear maps of vector spaces with fixed bases.
Each term of a sequence is a linear space of differentials of
metric values ascribed to the elements of a simplicial complex ---
a triangulation of a manifold. If the sequence turns out to be an acyclic
complex then one can construct a manifold invariant out of its
torsion. We demonstrate this first for three-dimensional manifolds,
and then we conduct the part, related to moves~$2\leftrightarrow 4$, of the
corresponding work for four-dimensional manifolds.
\end{abstract}

%%%%%%%%%%%%%
%\baselineskip1.85\baselineskip
%%%%%%%%%%%%%

\section{Introduction}

The present work is second in a series of papers whose aim is
to construct invariants of four-dimensional piecewise-linear manifolds
similar to the invariants of three-dimensional manifolds constructed
in papers \cite{3dim1} and~\cite{3dim2}. The first one was paper~\cite{I}
where we considered only rebuildings (Pachner moves) of type $3\to 3$,
that is replacements of three 4-simplices from a manifold triangulation
with three different simplices having the same common boundary. Below we
often use the Roman numeral~I when referring
to formulas, theorems, etc.\ from paper~\cite{I}
(and we call the paper~\cite{I} itself
paper~I). For instance, the expression for the invariant of moves $3\to 3$
is given in formula~(I.23).

It turns out that the further construction of our invariants involves
in a natural way the algebraic language of {\em acyclic complexes}.
Moreover, already in the three-dimensional case one of our constructions used in
papers \cite{3dim1} and~\cite{3dim2} has a natural interpretation
as the  calculation of the torsion of some properly constructed acyclic complexes.

Below in Section~\ref{sec b1.1} we study the three-dimensional case from this
new standpoint. One the one hand, this looks interesting in itself, and on the other
hand, the investigation of the three-dimensional case prepares the ground for
constructing and studying our algebraic complex for four-dimensional manifolds.
The general form of the complex and general plan of work in the four-dimensional
case are given in Section~\ref{sec b1.2}. Part of this plan is realized in this
paper, and the other part --- in the planned paper~III of this series.

In Section~\ref{sec b2.1} of the present paper we introduce {\em edge deviations\/}
that constitute one of the linear spaces in our complex. In Section~\ref{sec b2.2}
we partly justify the name ``complex'' by proving a theorem which states that
the image of one of the mappings lies in the kernel of the next mapping (we plan
to do this work in full, for the case of four dimensions, in paper~III). 
In Section~\ref{sec b2.3} we explain the relation between edge deviations and
moves $2\leftrightarrow 4$, i.e.\ replacements of a cluster of two 4-simplices
with a cluster of four ones and vice versa.

In the concluding Section~\ref{sec 2-4 disc} we sum up the results of paper~I and the present
paper and designate the aims of the next paper in this series.

\section{Revisiting the three-dimensional case: invariant~$I$ as
the torsion of an algebraic complex}
\label{sec b1.1}

We recall that in papers \cite{3dim1} and~\cite{3dim2} we were considering
orientable three-dimensional manifolds represented as simplicial
complexes (or ``pre-complexes'', where a simplex can enter several times
in the boundary of a simplex of a greater dimension). We
ascribed to every edge a Euclidean length, and to every 3-cell
(tetrahedron) --- a sign
$+$ or $-$ in such way that the following ``zero curvature'' condition
was satisfied: the algebraic sum of dihedral angles abutting at each edge
was zero modulo~$2\pi$. Here a dihedral angle is determined by the lengths of
edges of the relevant tetrahedron and is taken with the sign ascribed to that
tetrahedron. Then we considered infinitesimal variations of lengths
and edges and infinitesimal ``deficit angles'' at edges depending on them.
A matrix~$A$ that expressed these linear dependencies played a key r\^ole in
the resulting formula for manifold invariants.

It is natural to regard matrix~$A$ as a linear mapping of one vector space
with a fixed basis in it into another similar space. It turns out that $A$ makes up in a natural way, together
with four other linear maps of based linear spaces, an exact sequence,
or acyclic complex, that can be written out as follows:
\be
0
\leftarrow (\cdots)
\stackrel{B^{\rm T}}{\leftarrow} (d\omega) \stackrel{A}{\leftarrow} (dl)
\stackrel{B}{\leftarrow} (dx
\hbox{ and } dg) \leftarrow 0. 
\label{diagr} 
\ee 
Here $(dx)$ means all differentials of coordinates of vertices in the complex,
taken up to the motions of the Euclidean space (for example, in Section~3
of paper~\cite{3dim2} these were $\varphi$, $\sigma$, $s$ and $\alpha$);
$(dg)$ means all differentials of continuous parameters on which
a representation of the manifold's fundamental group in $E_3$ can depend
(in the case of an infinite fundamental group --- see explanations below);
$(dl)$ is the column of differentials of edge lengths;
$(d\omega)$ is the column of differentials of deficit angles. Recall that
the superscript $\rm T$ means the matrix transposing and that
$A=A^{\rm T}$. The left-hand half of the sequence~(\ref{diagr}) is obtained from
its right-hand half by the transposing of all matrices,
therefore, the exactness in terms ``$(\cdots)$'' and $(d\omega)$ follows from
the exactness in the respective terms in the right-hand half. We do not need
to determine the geometric meaning of the quantities denoted as
``$(\cdots)$'', --- it will follow from the sequence exactness
that they are simply elements of space $\Im B^{\rm T}$ --- the image of
map~$B^{\rm T}$.

We will need the following lemma.

\begin{lemma}
\label{lemma R3}
If edge lengths and tetrahedron signs in the complex are chosen so that
all deficit angles are zero then the vertices of its universal cover
can be put in a Euclidean space~$\mathbb R^3$ in such way that the lengths
will coincide with the distances between vertices, while tetrahedron signs --- with the orientations of
their images in~$\mathbb R^3$. This can be done in a unique way up to
the motions of space~$\mathbb R^3$.
\end{lemma}

{\it Proof}. Choose a tetrahedron in the universal cover --- we call it
the ``first'' tetrahedron --- and associate a Euclidean system of coordinates
with it. This system of coordinates can be naturally extended to any
adjoining tetrahedron, i.e.\ one having a common face with the first one
(recall that every {\em separate\/}
tetrahedron can always be put in~$\mathbb R^3$).
Proceeding further this way, we can introduce coordinates in any chosen
tetrahedron by using a sequence of tetrahedra starting at our first
tetrahedron, ending at our chosen ``last''
tetrahedron and such that any tetrahedron
in it adjoins the two neighbouring ones. We think of such a sequence
in geometric terms as determined by a broken line whose every straight segment
links the barycenters of two adjoining tetrahedra and whose beginning and end
are in the first and last tetrahedron respectively.

We are going to prove that if there are two such broken lines with the same
beginnings and ends then they yield the same system of coordinates in the last
tetrahedron. Assume that it is not so. Then we would have a closed path
with the beginning and end in the last tetrahedron and such that we obtain
a {\em changed\/} system of coordinates on going around that path. This path
is contractible (recall that we are in a universal cover). It can be contracted
into the barycenter of the last tetrahedron, and we can assume that in the process
of that contraction it passes through edges only (not vertices).

The fact that the deficit angle around every edge is zero implies that
Euclidean coordinates can be extended uniquely and without a contradiction from
any tetrahedron containing an edge to {\em all\/} tetrahedra
containing that edge. Thus, the monodromy matrix that describes the change
of coordinates corresponding to the way around our closed path does not
change when that path traverses an edge. Finally, we get a ``path'' of just
one point and two different systems of coordinates --- this contradiction
concludes the proof of the lemma.

\medskip

We will define the entries of sequence~(\ref{diagr}) and prove its exactness
with enough rigour only for manifolds~$M$ with a finite fundamental
group~$\pi_1(M)$, and we will only show what we can expect in the case
of infinite fundamental group on a simple example. Recall~\cite{3dim2}
that our invariant depends also on a representation $f\colon\; \pi_1(M)\to E_3$.
We start from the case considered in paper~\cite{3dim2}: $M=L(p,q)$ with
a nontrivial representation of group $\pi_1\bigl(L(p,q)\bigr)=\mathbb Z_p$
in $E_3$ by rotations around the $z$~axis.

On taking a universal cover, each vertex of the triangulation of $L(p,q)$
turns into $p$ copies of itself. We place those $p$ points, according to
paper~\cite{3dim2}, in the vertices of a regular $p$-gon in such way that they
turn into one another under rotations through multiples of $2\pi/p$
around the $z$~axis. In contrast with paper~\cite{3dim2}, we are now considering,
however, an {\em arbitrary\/} triangulation --- having any number of vertices.
Let there be $m$ vertices, then their positions with respect to each other
are determined, first, by $m$ distances $\rho_1,\ldots,\rho_m$ between them
and the $z$~axis and, second, by $(m-1)$ differences for each of two remaining
cylindrical coordinates, e.g., $(\varphi_2-\varphi_1),\ldots, (\varphi_m-\varphi_1)$
and $(z_2-z_1),\ldots, (z_m-z_1)$. The experience of paper~\cite{3dim2} shows
that the differentials $dx$ must be chosen as (cf.\
formula~\cite[(3.16)]{3dim2})
\be
(dx)=\bigl(\rho_1\, d\rho_1,\ldots, \rho_m\, d\rho_m;
d(\varphi_2-\varphi_1),\ldots, d(\varphi_m-\varphi_1); d(z_2-z_1),\ldots,
d(z_m-z_1)\bigr).
\label{rho phi z}
\ee

As for the differentials $(dg)$, there are none of them in this case:
a representation of a finite group cannot be deformed continuously
in a non-equivalent one.

Map~$B$ determines the deformations of all edge lengths in the complex for
given~$(dx)$. Its injectivity for generic $\rho$, $\varphi$ and $z$ is
obvious: the correspondence between $\rho_i$, $\varphi_i-\varphi_1$,
$z_i-z_1$, on the one hand, and the edge lengths in the complex
yielding zero deficit angles, is locally one-to-one (this follows easily from
Lemma~\ref{lemma R3}), and thus the Jacobian determinant of map~$B$ is almost everywhere
nonzero. The injectivity of~$B$ means the exactness of
sequence~(\ref{diagr}) in the term~$(dx)$.

So, it remains to prove the exactness in the term~$(dl)$. The fact that
$\Im B \subset \Ker A$ is obvious: if length differentials are generated by
translations of vertices within the Euclidean space~$\mathbb R^3$, then
the deficit angles are zero. Let us prove that $\Ker A \subset \Im B$. This
means that if infinitesimal deficit angles are zero then our complex --- the
universal cover of a triangulation of~$L(p,q)$ --- can be put in the
Euclidean space~$\mathbb R^3$. This is, however, nothing else but an
infinitesimal version of the same Lemma~\ref{lemma R3} (with the same
proof).

Thus, the exactness of sequence~(\ref{diagr}) for $L(p,q)$ is proven.
We would like to generalize this result at once for all (connected
closed orientable) manifolds~$M$ with a {\em finite\/} fundamental
group $\pi_1(M)$. A single important difference from the case of $L(p,q)$ can
consist in the choice of variables~$(dx)$. As one can see from
papers~\cite{3dim1,3dim2}, the variables~$(dx)$ must possess the following
property: on adding a new vertex~$E$ (as a result of a move $1\to 4$) the
form $\bigwedge dx$ must get multiplied by $dx_E \wedge dy_E \wedge dz_E$ ---
the form of Euclidean volume. Recall that our invariant corresponds to a
{\em pair\/} $(M,f)$, where $f\colon\; \pi_1(M)\to E_3$. Three situations are
possible depending on the {\em image\/} $\Im f$ of the group $\pi_1(M)$
in~$E_3$:

1)~$\Im f=\{e\}$. This case has been studied in paper~\cite{3dim1}.

2)~$\Im f$ contains only rotations around one axis. In this case, the
same reasonings as presented above for $L(p,q)$ are valid.
In particular, $(dx)$ must be chosen
according to formula~(\ref{rho phi z}).

3)~$\Im f$ contains rotations around two or more (nonparallel) axes.
The presence of such axes fixes the system of coordinates in~$\mathbb R^3$, and
$(dx)$ must simply consist of all differentials of Euclidean coordinates
of vertices in the complex:
\be
(dx)=(dx_1,dy_1,dz_1,\ldots,dx_m,dy_m,dz_m).
\label{dx sluchaj 3}
\ee

Having presented an acyclic complex (\ref{diagr}) for a manifold~$M$ with
a finite group $\pi_1(M)$, we now consider its {\em torsion\/}. As in
papers \cite{3dim1,3dim2}, we choose a maximal subset $\cal C$ in the set of
edges of $M$'s triangulation such that the square matrix $A|_{\cal C}$
(obtained from $A$ by removing the rows and columns corresponding to
the edges from~$\ol {\cal C}$ --- the complement of~$\cal C$) is nondegenerate.
We also need a matrix~$B|_{\ol{\cal C}}$ obtained by removing from
matrix~$B$ the {\em rows} corresponding to the edges from~$\cal C$.
According to usual formulae~\cite{kruchenie}, the torsion~$\tau$ is
\be
\tau=\frac{\left( \det B|_{\ol {\cal C}} \right)^2}{\det A|_{\cal C}}\,.
\label{tau 3d}
\ee
As is well-known, $\tau$ does not depend on the choice of~$\cal C$.

We can write in the style of papers \cite{3dim1,3dim2}:
$$
\det B|_{\ol {\cal C}} = \frac{\displaystyle \bigwedge_{\ol{\cal C}} dl}
{\displaystyle \bigwedge dx}\,.
$$
After this, the translation of those works into our new language presents
no difficulty. We formulate the result as the following Theorem.

\begin{theorem}
\label{th inv 3d}
The value
\be
\frac{\displaystyle \tau \prod_{\rm over\; all\; edges} l^2}{\displaystyle
\prod_{\rm over\; all\; tetrahedra} 6V}
\label{inv 3d}
\ee
is an invariant of a three-dimensional connected closed orientable
manifold with a finite fundamental group.
\end{theorem}

Clearly, (\ref{inv 3d}) is the {\em squared\/} invariant $I$ from \cite{3dim1}
or $I_k$ from \cite{3dim2} in the cases considered in those papers.

\medskip

{\it Proof\/} of Theorem \ref{th inv 3d} consists in tracing what happens
with expression (\ref{inv 3d}) under moves $2\leftrightarrow 3$ and
$1\leftrightarrow 4$, and this is done in the very same way as in
papers \cite{3dim1} and \cite{3dim2}. The theorem is proven.

\medskip

We show now what we can expect in the case of an infinite $\pi_1(M)$ on the
example of a manifold~$M$ with $\pi_1(M)=\mathbb Z$ (for instance,
$M=S^1\times S^2$). Generically, the generating element of $\pi_1(M)$ is
sent by homomorphism~$f$ in a rotation around some axis --- let that be
the $z$~axis --- through an angle~$\alpha$ and a translation along the
same axis through a distance~$a$. Arguments similar to those given above
show that we will get an exact sequence~(\ref{diagr}) if we add, in the
capacity of $(dg)$, the differentials $d\alpha$ and $da$ to the same $(dx)$
as in formula~(\ref{rho phi z}). After this, the torsion and invariant are
calculated according to the old formulas (\ref{tau 3d}) and~(\ref{inv 3d}).
We are not giving the details here
and just put forward a conjecture that one can choose
$(dx)$ and $(dg)$ for any $M$ so that (\ref{diagr}) will be an exact sequence.

\section{Four-dimensional case: the algebraic complex
in its general form}
\label{sec b1.2}

We will now present a sequence of linear maps of based linear spaces which
is our candidate for the r\^ole of acyclic complex, similar
to~(\ref{diagr}), in the four-dimensional case. To be exact, we write down
two ``conjugate'' sequences: the matrices defining the maps in one of them
can be obtained by matrix transposing from the matrices in the other sequence.
This enables us, while investigating the exactness in various terms
of the sequences, to use that sequence which is more suitable for the
given case. Moreover, this will allow us to bypass the question about the
geometrical meaning of some linear spaces entering in a sequence if the
meaning of the corresponding fragment of the other sequence is clear and that
fragment lends itself to investigation. We will denote such linear spaces
by the marks of omission~$(\cdots)$. So, we are using this notation for {\em different\/}
linear spaces which are not (yet) involved in our reasonings.

The definitions of some linear spaces entering in the first of our
sequences will be given in the next sections of the present paper and
in the forthcoming papers of this series. Nevertheless, we believe it
reasonable to write down the sequences right now, because they will
determine the plan of our further work and the interrelations of its
individual parts.

So, here are our sequences:
\begin{eqnarray}
0\leftarrow (\cdots) \leftarrow (d\Omega_a) \stackrel{(\partial
\Omega_a / \partial S_i)}{\longleftarrow} (dS_i) \leftarrow (d\vec v_a)
\leftarrow (d\sigma) \leftarrow 0\,,
\label{diag4} \\
0 \rightarrow (dx \hbox{ and } dg) \rightarrow (dL_a) \stackrel{(\partial
\omega_i / \partial L_a)}{\longrightarrow} (d\omega_i) \rightarrow
(\cdots) \rightarrow (\cdots) \rightarrow 0\,.
\label{diag4T}
\end{eqnarray}
We start the explanations with the sequence (\ref{diag4T}). Its left-hand
side to $(d\omega_i)$ inclusive is a direct analogue of the
{\em right-hand\/} side of sequence~(\ref{diagr}). We only recall that
$L_a$ is the {\em squared\/} length of the $a$th edge, while $\omega_i$ is
the deficit angle around the $i$th {\em two-dimensional face}.

The link between (\ref{diag4T}) and (\ref{diag4}) is provided by formula~(I.16)
which states the mutual conjugacy of matrices $(\partial \Omega_a / \partial
S_i)$ and~$(\partial \omega_i / \partial L_a)$. Recall that $S_i$ is the area
of $i$th two-dimensional face, while $\Omega_i$ is the deficit angle around
the $a$th edge defined according to formula~(I.13). The values $d\vec v_a$ ---
we call them {\em edge deviations} --- are introduced in Section~\ref{sec b2.1},
together with the mapping $(d\vec v_a)\to (dS_i)$. As for the values
$d\sigma$ --- {\em vertex deviations} --- we leave their definition to
the future paper~III of this series.

It seems that there exist enough interesting manifolds for which our
sequences (\ref{diag4}) and~(\ref{diag4T}) are acyclic complexes and hence
one can hope to obtain four-dimensional manifold invariants out of their torsion
in analogy with formula~(\ref{inv 3d}). In the present paper, we concentrate on
the information extractable from fragment $(d\Omega_a)\leftarrow (dS_i)\leftarrow (d\vec v_a)$
of sequence~(\ref{diag4}). This fragment turns out to be responsible, in a sense,
for moves~$2\leftrightarrow 4$.

\section{Edge deviations and area differentials determined by them}
\label{sec b2.1}

In this Section we introduce the notion of edge deviation and explain how
these deviations give rise to differentials of two-dimensional face areas,
or, in other words, how the matrix~$(\partial S_i / \partial \vec v_a)$
is constructed that gives the map $(d\vec v_a)\to (dS_i)$ in
sequence~(\ref{diag4}).

Let $AB$ be some edge of our complex. Recall that in Section~3 of paper~I we have
mapped all cells of the {\em universal cover\/} of the complex
in a Euclidean space~$\mathbb R^4$. Using the pull-back of edge~$AB$ and adjoining
4-simplices up to the universal cover, we can assume that this edge and all neighbouring
elements of the complex are put in~$\mathbb R^4$. We call the {\em deviation of edge\/}~$AB$
any infinitesimal vector $d\vec v_{AB}$ orthogonal to~$AB$.
We imagine such vector as having its beginning in the {\em middle of edge~$AB$}.

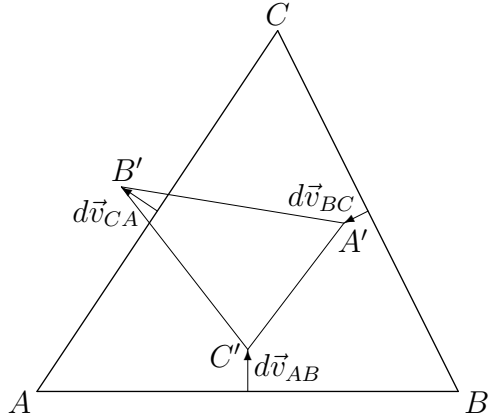
\begin{figure}
\begin{center}
\unitlength=0.8mm
\special{em:linewidth 0.5pt}
\begin{picture}(76.00,66.00)
\emline{5.00}{5.00}{1}{45.00}{65.00}{2}
\emline{45.00}{65.00}{3}{75.00}{5.00}{4}
\emline{75.00}{5.00}{5}{5.00}{5.00}{6}
\special{em:linewidth 0.2pt}
\thinlines
\emline{40.00}{5.00}{7}{40.00}{12.00}{8}
\emline{60.00}{35.00}{9}{56.00}{33.00}{10}
\emline{25.00}{35.00}{11}{19.00}{39.00}{12}
\emline{40.00}{12.00}{13}{56.00}{33.00}{14}
\emline{56.00}{33.00}{15}{19.00}{39.00}{16}
\emline{19.00}{39.00}{17}{40.00}{12.00}{18}
\put(56.00,33.00){\vector(-2,-1){0.00}}
\put(19.00,39.00){\vector(-3,2){0.00}}
\put(40.00,12.00){\vector(0,1){0.00}}
\put(4.00,5.00){\makebox(0,0)[rt]{$A$}}
\put(76.00,5.00){\makebox(0,0)[lt]{$B$}}
\put(45.00,66.00){\makebox(0,0)[cb]{$C$}}
\put(58.00,35.00){\makebox(0,0)[rb]{$d\vec v_{BC}$}}
\put(41.00,9.00){\makebox(0,0)[lc]{$d\vec v_{AB}$}}
\put(22.00,37.00){\makebox(0,0)[rt]{$d\vec v_{CA}$}}
\put(20.00,40.00){\makebox(0,0)[cb]{$B'$}}
\put(55.00,32.00){\makebox(0,0)[lt]{$A'$}}
\put(39.00,13.00){\makebox(0,0)[rt]{$C'$}}
\end{picture}
\end{center}
\caption{Illustration to formula~(\protect\ref{dSABC})}
\label{fig devABC}
\end{figure}

The area of any Euclidean triangle is four times area of the triangle
with vertices in the middles of its sides. If edge deviations in
triangle~$ABC$ are $d\vec v_{AB}$, $d\vec v_{AC}$ and $d\vec v_{BC}$ then
we set by definition that the differential $dS_{ABC}$ generated by them is
\be
dS_{ABC}=4S_{A'B'C'}-S_{ABC}
\label{dSABC}
\ee
(Figure~\ref{fig devABC}).

Now let there be a 4-simplex $ABCDE$ lying in Euclidean space $\mathbb R^4$
with coordinates $(x,y,z,t)$. Assume that vertices $A$ and~$B$ lie on the
$t$~axis. Consider a situation where only edge $AB$ has a nonzero deviation
$$
d\vec v_{AB} = d\vec v = (dv_x, dv_y, dv_z, 0).
$$
We are going to derive a formula that relates
$dv_x \wedge dv_y \wedge dv_z$ to
$dS_{ABC} \wedge dS_{ABD} \wedge dS_{ABE}$.

We project simplex $ABCDE$ along the $t$~axis onto a {\em three-dimensional\/}
space and draw the result in thick lines in Figure~\ref{fig deviacii reber},
\begin{figure}
\begin{center}
\unitlength=0.80mm
\special{em:linewidth 0.5pt}
\begin{picture}(76.00,78.00)
\emline{35.00}{5.00}{1}{75.00}{25.00}{2}
\emline{75.00}{25.00}{3}{45.00}{75.00}{4}
\emline{45.00}{75.00}{5}{5.00}{25.00}{6}
\emline{5.00}{25.00}{7}{35.00}{5.00}{8}
\emline{35.00}{5.00}{9}{45.00}{75.00}{10}
\emline{5.00}{25.00}{11}{11.00}{25.00}{12}
\emline{16.00}{25.00}{13}{22.00}{25.00}{14}
\emline{27.00}{25.00}{15}{33.00}{25.00}{16}
\emline{38.00}{25.00}{17}{44.00}{25.00}{18}
\emline{49.00}{25.00}{19}{55.00}{25.00}{20}
\emline{60.00}{25.00}{21}{66.00}{25.00}{22}
\emline{71.00}{25.00}{23}{75.00}{25.00}{24}
\special{em:linewidth 0.2pt}
\thinlines
\put(56.00,72.00){\vector(3,-1){00.00}}
\emline{25.00}{50.00}{25}{40.00}{40.00}{26}
\emline{40.00}{40.00}{27}{60.00}{50.00}{28}
\emline{60.00}{50.00}{29}{56.00}{72.00}{30}
\emline{56.00}{72.00}{31}{25.00}{50.00}{32}
\emline{40.00}{40.00}{33}{56.00}{72.00}{34}
\emline{25.00}{50.00}{35}{27.00}{50.00}{36}
\emline{30.00}{50.00}{37}{34.00}{50.00}{38}
\emline{37.00}{50.00}{39}{41.00}{50.00}{40}
\emline{44.00}{50.00}{41}{48.00}{50.00}{42}
\emline{51.00}{50.00}{43}{55.00}{50.00}{44}
\emline{58.00}{50.00}{45}{60.00}{50.00}{46}
\put(44.00,74.00){\makebox(0,0)[rb]{$A,B$}}
\put(49.00,75.00){\makebox(0,0)[lb]{$d\vec v$}}
\put(57.00,72.00){\makebox(0,0)[lc]{$K$}}
\put(61.00,50.00){\makebox(0,0)[lb]{$E_{1/2}$}}
\put(41.00,40.00){\makebox(0,0)[lt]{$D_{1/2}$}}
\put(24.00,50.00){\makebox(0,0)[rb]{$C_{1/2}$}}
\put(4.00,25.00){\makebox(0,0)[rc]{$C$}}
\put(35.00,4.00){\makebox(0,0)[ct]{$D$}}
\put(76.00,25.00){\makebox(0,0)[lc]{$E$}}
\emline{45.00}{75.00}{47}{56.00}{72.00}{48}
\end{picture}
\end{center}
\caption{Three-dimensional projection of 4-simplex $ABCDE$ and the deviation
of edge~$AB$}
\label{fig deviacii reber}
\end{figure}
denoting somewhat loosely the {\em projections\/} of points $A,\ldots,E$ by
the same letters. The other points in Figure~\ref{fig deviacii reber} are:
$C_{1/2}$ is the common projection of middle points of edges $AC$ and~$BC$;
similarly $D_{1/2}$ is the projection of middle points of edges $AD$ and~$BD$,
and $E_{1/2}$ of edges $AE$ and~$BE$;
$K$ is the projection of the point distant from the middle of $AB$
by the deviation~$d\vec v$.

Forgetting for a moment about the four-dimensional origin of
Figure~\ref{fig deviacii reber}, we can write down, using an easy
trigonometry (cf.\ formulas (31) and (32) of paper~\cite{3dim1}):
$$
|dv_x \wedge dv_y \wedge dv_z| = \left| \frac{
l_{AC_{1/2}}\, dl_{KC_{1/2}} \wedge
l_{AD_{1/2}}\, dl_{KD_{1/2}} \wedge
l_{AE_{1/2}}\, dl_{KE_{1/2}} }
{6V_{AC_{1/2}D_{1/2}E_{1/2}}} \right|.
$$
Here, of course, $dl_{KC_{1/2}}=l_{KC_{1/2}}-l_{AC_{1/2}}$ and so~on;
the letter~$l$ itself denotes the length of the corresponding edge lying in the
{\em three-dimensional\/} space, while $V$ is the three-dimensional volume.

Now we go back into the fourth dimension. Denote the length of edge~$AB$ as~$h$.
It follows from the orthogonality of this edge to the three-dimensional projection
depicted in Figure~\ref{fig deviacii reber} that
$$
| 6V_{AC_{1/2}D_{1/2}E_{1/2}}\cdot h | = | 3V_{ABCDE} |
$$
(as usual, we use the absolute value signs in order not to care about
the orientation),
\be
l_{AC_{1/2}}\cdot h = S_{ABC},\quad dl_{KC_{1/2}}\cdot h = dS_{ABC},
\label{lAC1/2*h}
\ee
and two more pairs of formulas of type (\ref{lAC1/2*h}) are got by changes
$C\to D$ and~$C\to E$.

All this together gives
\be
L_{AB}^{5/2} \left|dv_x \wedge dv_y \wedge dv_z \right| =
\left| \frac { d\left(S_{ABC}^2\right) \wedge d\left(S_{ABD}^2\right) \wedge
d\left(S_{ABE}^2\right) } {24\, V_{ABCDE}} \right|,
\label{devAB}
\ee
where we have returned to the notation $L_{AB}=h^2$ of paper~I for a squared
edge length.

\section{Edge deviations yield zero deficit angles~$d\Omega$}
\label{sec b2.2}

The following Theorem gives ``half'' of the exactness of sequence
(\ref{diag4}) in the term~$(dS_i)$.

\begin{theorem}
\label{th 1/2 v Si}
In sequence (\ref{diag4}), $\Im (\partial S_i / \partial \vec v_a)
\subset \Ker (\partial \Omega_a / \partial S_i )$.
\end{theorem}

{\it Proof} constitutes all the remaining part of this Subsection.

We represent the totality of all differentials $dS_i$ of two-dimensional
face areas in the complex as a column vector~$(dS_i)$. We consider first
those $dS_i$ that arise from deviations of {\em one\/} edge~$AB$ common
for {\em exactly four\/} 4-simplices $ABCDE=\hat F$, $ABCFD=-\hat E$, $ABCEF=\hat D$
and~$ABDFE=-\hat C$ (written in such form, these simplices have a {\em
consistent orientation\/}). As $AB$ is common also for exactly four
two-dimensional faces, vector~$(dS_i)$ can have in our case only
four nonzero components. They turn out to obey the following linear
dependence:
\be
V_{-\hat C}\, d(S_{ABC}^2) + V_{\hat D}\, d(S_{ABD}^2) +
V_{-\hat E}\, d(S_{ABE}^2) + V_{\hat F}\, d(S_{ABF}^2) = 0\,.
\label{VdS}
\ee

Indeed, we see from formula (\ref{devAB}) that the expression
\be
| d(S_{ABC}^2) \wedge d(S_{ABD}^2) \wedge d(S_{ABE}^2) / V_{\hat F} |
\label{VdS*}
\ee
is invariant with respect to permutations of letters $C,D,E,F$. One can check
that~(\ref{VdS}) is exactly the linear dependence that ensures
this invariance. To be exact, the summands in (\ref{VdS}) are determined from
the invariance of (\ref{VdS*}) up to their signs, but one can easily fix
those signs by considering, e.g., such limit cases when one of points
$D$, $E$ or $F$ comes very closely to point~$C$: if $D\to C$, then the
differentials $dS_{ABC}$ and $dS_{ABD}$ must be identical.

We prove that~(\ref{VdS}) is equivalent to the condition $d\Omega_{AB}=0$,
where $d\Omega_{AB}$ --- the deficit angle around edge~$AB$ --- is considered
as a function of $dS_{ABC}$, $dS_{ABD}$, $dS_{ABE}$ and~$dS_{ABF}$. Indeed,
it follows from the mutual conjugacy of matrices
$(\partial \Omega_a / \partial S_i)$ and
$(\partial \omega_i / \partial L_a)$ (Theorem~I.2) that a vector
$(dS_i)$ that yields zero deficit angles around edges must be {\em
orthogonal\/} to the image of matrix $(\partial \omega_i / \partial L_a)$,
i.e.\ to all its columns. Considering the column corresponding to edge $a=AB$,
and applying formula (I.5) and formulas obtained from it by changes
$C\leftrightarrow D$, $C\leftrightarrow E$ and $C\leftrightarrow F$, we get:
\be
\pa{\omega_{ABC}}{L_{AB}} : \pa{\omega_{ABD}}{L_{AB}} :
\pa{\omega_{ABE}}{L_{AB}} : \pa{\omega_{ABF}}{L_{AB}} =
V_{-\hat C}S_{ABC} : V_{\hat D}S_{ABD} :
V_{-\hat E}S_{ABE} : V_{\hat F}S_{ABF}\,.
\label{6:}
\ee
This leads at once to linear dependence (\ref{VdS}), and we have only to note
that {\em only one\/} dependence exists between our four $dS$.

The fact that equality~(\ref{VdS}), on the one hand, holds for area differentials
generated by deviations of edge~$AB$ and, on the other hand, is equivalent to
the condition $d\Omega_{AB}=0$, shows that we have proved the following lemma.

\begin{lemma}
\label{lemma AB}
If an edge $AB$ is common for exactly four 4-simplices, then the infinitesimal
changes of areas arising from deviations $d\vec v_{AB}$ give rise, in their
turn, to a zero deficit angle $d\Omega_{AB}=0$.
\end{lemma}

We want to prove that deviations of {\em all\/} edges in the complex give rise
to such $dS_i$ that {\em all\/} $d\Omega_a=0$. Let us prove this first
for a complex representing a triangulation of sphere~$S^4$ consisting of
six vertices $A,\ldots,F$ and six 4-simplices $\hat A,\ldots,\hat F$.

In such complex, first, {\em every\/} edge is common for exactly four
4-simplices. Second, all $d\Omega_a$ are proportional (the rank of matrix
$(\partial \Omega_a / \partial S_i)$ is~1). Thus, Lemma~\ref{lemma AB} gives
at once the needed result.

If we now have an arbitrary complex that contains an edge~$AB$ common for
exactly four 4-simplices, then, again, deviations of {\em all\/} edges yield
$d\Omega_{AB}=0$, because $d\Omega_{AB}$ can depend, in principle, only on
the deviations of the edges belonging to the four ``nearby'' simplices,
and this dependence is the same as in the previous case of the
triangulation of~$S^4$.

It remains to consider only one more case, with more than four 4-simplices
situated around the edge~$AB$. We can always assume that our manifold
triangulation was obtained by means of Pachner moves, i.e.\ rebuildings
$3\to 3$, $2\leftrightarrow 4$ and $1\leftrightarrow 5$, from another
triangulation where {\em there was no\/} edge~$AB$.
That edge appeared as a result of a move $2\to 4$ or $1\to 5$, which means
that it was ``at first'' surrounded by four 4-simplices, so that all
the derivatives $(\partial \Omega_{AB} / \partial \vec v_a)$ were zero.
Then, more rebuildings of various types could be performed in a neighbourhood
of edge~$AB$, but any of them may be interpreted in the following way:
a cluster of 4-simplices in our simplicial complex was replaced with
another cluster in such way that they both formed together the
triangulation of sphere~$S^4$ consisting of
six 4-simplices (as described above).
One may say that six new 4-simplices are added to the complex, and then
some of them ``cancel out'' with those already present, because
their orientations are opposite. Together with six added simplices
(whose orientations are {\em consistent between themselves\/} for any
rebuilding), the deficit angles pertaining
to~$S^4$ are added to the corresponding
deficit angles in the complex. Since all derivatives
$\partial \Omega_{AB} / \partial \vec v_a$ are zero in both $S^4$ and the
``old'' complex, they remain zero after the rebuilding as well.

Thus, Theorem~\ref{th 1/2 v Si} is proven.

\section{Moves $2\leftrightarrow 4$ and the invariant of moves $3\to 3$}
\label{sec b2.3}

In paper I, we have constructed the invariant (I.23) of moves $3\to 3$ {\em only\/}
in the form of the product of two differential forms (acting on exterior
powers of two different linear spaces). From the standpoint of our sequences (\ref{diag4})
and~(\ref{diag4T}), we were studying only the map $(\partial \Omega_a / \partial S_i)$
and its conjugate $(\partial \omega_i / \partial L_a)$. Expression~(I.23)
cannot be invariant with respect to moves $2\leftrightarrow 4$ because they change
the degree of form~$\bigwedge dS$ --- three differentials are added or taken away. Still,
we are going to show that, assuming the exactness of sequence~(\ref{diag4})
in the term~$(dS_i)$, we can describe the behaviour of expression~(I.23)
under moves $2\leftrightarrow 4$ in quite simple terms.

For concreteness, we speak below about the replacement of two 4-simplices with
four ones, i.e.\ a move~$2\to 4$. All our reasonings and formulas admit a reversion
so that they can also describe a move~$4\to 2$.

So, let two adjoining 4-simplices $ACDEF=\hat B$ and $BCDFE=-\hat A$ be replaced
by four of them: $ABCDE=\hat F$, $ABCFD=-\hat E$, $ABCEF=\hat D$
and~$ABDFE=-\hat C$. This rebuilding provides the complex with a new
edge~$AB$ and four new two-dimensional faces containing that edge.
Arguments similar to those in Section~3 of paper~\cite{3dim1} show that
the rank of matrix $(\partial \omega_i / \partial L_a)$ increases thus by~1
and, if edge~$AB$ is added to the subset~$\cal C$ of the set of all edges
(see~I, Section~5) and, for instance, face $ABF$ is added to the subset
$\cal D$ of the set of all faces (ibidem), then the determinant of matrix ${\cal B} =
\left.\vphantom{|}\right._{\cal D}\!
(\partial \omega_i / \, \partial L_a )
\!\!\left.\vphantom{|}\right._{\cal C}$
gets multiplied under our rebuilding by $\pm \partial \omega_{ABF} / \partial L_{AB}$:
\be
\det {\cal B} \mathrel{\mathop{\longrightarrow}\limits_{2\to 4}}
\pm \pa{\omega_{ABF}}{L_{AB}} \det {\cal B}\,.
\label{detB->}
\ee

As for the partial derivative in~(\ref{detB->}), it can be obtained by
comparing formula~(I.5) for a similar derivative $\partial \omega_{ABC} / \partial L_{AB}$
with formula~(\ref{6:}):
\be
\pa{\omega_{ABF}}{L_{AB}} = \frac{S_{ABF}}{24}
\frac{V_{\hat A}\, V_{\hat B}}{V_{\hat C}\,V_{\hat D}\,V_{\hat E}}\,.
\label{ABF/AB}
\ee

Using formulas (\ref{detB->}) and (\ref{ABF/AB}), we find that invariant~(I.23)
gets multiplied under a move $2\to 4$ by
\be
\pm 3\, \frac{d(S_{ABC}^2) \wedge d(S_{ABD}^2) \wedge d(S_{ABE}^2)}{V_{\hat F}}
\label{SSS/V}
\ee
(of course, the {\em exterior\/} product of new area differentials by those
already existing is taken).

Recall now our assumption about the exactness of sequence~(\ref{diag4})
in the term~$(dS_i)$. We will need the inclusion inverse to that proven in
Theorem~\ref{th 1/2 v Si}, namely,
\be
\Ker \left(\pa{\Omega_a}{S_i}\right) \subset \Im \left(\pa{S_i}{\vec v_a}
\right).
\label{vkl obr}
\ee
It means that any column vector~$(dS_i)$ yielding zero $d\Omega_a$
can be obtained from deviations~$d\vec v_a$. Therefore, the whole
exterior product $\bigwedge dS$ in invariant~(I.23) can be expressed
through edge deviations.

Inclusion (\ref{vkl obr}) is conserved under a move $2\to 4$: one can always
choose three components of the new deviation~$d\vec v_{AB}$ so as to get
any values of three new independent (under the conditions $d\Omega_a=0$)
area differentials. The new $dS$ (entering in formula~(\ref{SSS/V})) depend
not only on $d\vec v_{AB}$ but also on other deviations, but the dependence
on those latter plays no r\^ole when taking the exterior
product~$\bigwedge dS$ in~(I.23), because $\bigwedge dS$ has always
the maximal degree, coinciding with the maximal number of linearly independent
$dS$ that can be obtained from edge deviations. Thus, the exterior product
(\ref{SSS/V}) can be re-written with the help of formula~(\ref{devAB}). Hence,
we get the following Theorem.

\begin{theorem}
\label{th 3-3 2-4}
Assuming that inclusion (\ref{vkl obr}) holds, invariant (I.23)
of moves $3\to 3$ gets multiplied under a move $2\to 4$ by
\be
\pm 72\, L_{AB}^{5/2} \, dv_x \wedge dv_y \wedge dv_z\,,
\label{SSS/V*}
\ee
where $AB$ is the added edge, while $dv_x$, $dv_y$ and~$dv_z$ are the components
of its deviation.
\end{theorem}

\section{Discussion}
\label{sec 2-4 disc}

Theorem~\ref{th 3-3 2-4} shows that the invariant of moves $3\to 3$
introduced in paper~\cite{I} undergoes in fact only very simple changes
under moves $2\leftrightarrow 4$ as well, at least if we assume
the inclusion~(\ref{vkl obr}). One can check that this inclusion is
satisfied for the sphere~$S^4$. Theorem~\ref{th 3-3 2-4} shows also that
it does make sense to consider at least the part of our sequence~(\ref{diag4})
from $(d\Omega_a)$ to~$(d\vec v_a)$ inclusive.

In the next, third paper of this series we plan to construct a linear space of
``vertex deviations'' $d\sigma$ together with the linear map from it to the space of
edge deviations~$d\vec v$. This map will be intimately connected with
moves~$1\leftrightarrow 5$. Thus, we will be able to take into account
{\em all\/} Pachner moves and construct the full invariant of a
four-dimensional piecewise-linear manifold with the help of the
torsion of sequence~(\ref{diag4}). This requires that the latter should be an
acyclic complex --- we will prove this fact, in particular, for the
sphere~$S^4$. We will show also that the exactness property itself is conserved
under Pachner moves.

Note that the investigation of sequence~(\ref{diag4}) looks very interesting
even if it fails to be exact for some manifold. According to some
arguments (suggested by Wall's stabilization theorem), this may
happen for the product $S^2 \times S^2$ of two-dimensional spheres.

In passing, we have translated some constructions in papers~\cite{3dim1,3dim2}
into the language of the torsion of acyclic complexes, and thus clarified
their algebraic nature. Schematically, the main characters in our play are,
first, solutions to the pentagon equation or its generalizations
(in dimensions $>3$) and, second, acyclic complexes. In this connection, it
would be very interesting to construct similar complexes corresponding to
the ``$SL(2)$ solution'' of the pentagon equation presented in
paper~\cite{SL_2} and its analogues for higher dimensions.

We conclude this paper with a question which looks the most intriguing:
what are the quantum objects whose limit yields our constructions,
with their so clearly quasiclassical appearance.

\medskip

{\bf Acknowledgements. }\vadjust{\nobreak}The work has been
performed with a partial financial support from Russian
Foundation for Basic Research under Grant no.~01-01-00059.


\begin{thebibliography}{99}

\bibitem{I}
{\it I. G. Korepanov}.
Euclidean 4-simplices and invariants of four-dimensional manifolds:
I.~Moves $3\to 3$.
Theor. Math. Phys., Volume~131, no.~3 (2002), 765--774.

\bibitem{3dim1}
{\it I. G. Korepanov}.
Invariants of PL manifolds from metrized simplicial
complexes. Three-dimensional case.
J. Nonlin. Math. Phys., Volume 8, number 2 (2001), 196--210.

\bibitem{3dim2}
{\it I. G. Korepanov, E. V. Martyushev}.
Distinguishing three-dimensional lens spaces $L(7,1)$ and $L(7,2)$
by means of classical pentagon equation.
J. Nonlin. Math. Phys., Volume 9, number 1 (2002), 86--98.

\bibitem{kruchenie}
{\it J. W. Milnor}. Whitehead torsion.
Bull. Amer. Math. Soc. Volume~72 (1966), 358--426.

\bibitem{SL_2}
{\it I. G. Korepanov, E. V. Martyushev}.
A Classical Solution of the Pentagon Equation Related to the Group $SL(2)$.
Theor. Math. Phys., Volume~129 (2001), 1320--1324.

\end{thebibliography}
\end{document}